\documentclass[11pt,reqno]{amsart}
\usepackage[T1]{fontenc}\usepackage{lmodern}
\usepackage[utf8]{inputenc}
\usepackage{amsmath,amssymb,mathrsfs}
\usepackage[left=0.95in,right=0.95in,top=0.95in]{geometry}
\usepackage{microtype}
\usepackage{xcolor}
\definecolor{linkcyan}{RGB}{0,112,128}
\usepackage[colorlinks=true,linkcolor=linkcyan,citecolor=linkcyan,urlcolor=linkcyan]{hyperref}
\hypersetup{
 pdftitle={There exist no embedded index six minimal surfaces in the round three-sphere},
 pdfauthor={Alessandro Carlotto}
}
\newtheorem{theorem}{Theorem}[section]
\newtheorem{proposition}[theorem]{Proposition}
\newtheorem{lemma}[theorem]{Lemma}
\newtheorem{corollary}[theorem]{Corollary}
\newtheorem{conjecture}[theorem]{Conjecture}
\theoremstyle{remark}
\newtheorem{remark}[theorem]{Remark}
\DeclareMathOperator{\ind}{ind}

\DeclareMathOperator{\diver}{div}
\DeclareMathOperator{\tr}{tr}
\DeclareMathOperator{\Sym}{Sym}
\newcommand{\R}{\mathbb R}
\newcommand{\HH}{\mathcal H}
\newcommand{\dd}{\,dA_h}
\title[No embedded index six minimal surfaces in the three-sphere]
{There exist no embedded index six minimal surfaces in the round three-sphere}
\author{Alessandro Carlotto}
\address{Universit\`a di Trento, Dipartimento di Matematica,
via Sommarive 14, 38123 Trento, Italy}
\email{alessandro.carlotto@unitn.it}
\date{}
\begin{document}
\begin{abstract}
We improve Savo's lower bound for the Morse index of closed
orientable minimal immersions in the round
three-sphere, in terms of the genus and nullity. As a consequence, we are able to show that
no closed embedded minimal surface has Morse index equal to six. 
\end{abstract}
\maketitle

\section{Introduction}\label{sec:introduction}

Since Lawson's pioneering construction in 1970 \cite{Lawson}, the landscape of minimal surfaces in the round three-sphere has grown enormously, and several classes of new examples have been obtained by variational, desingularization and spectral methods. Such a scenario naturally calls for classification results, and the present article presents a contribution in that direction. Here we shall be concerned with the study of minimal surfaces based on the value of their Morse index $I$, that is, the maximal dimension of a space of normal deformations on which
the second variation of area is negative definite.

For a closed, connected, orientable minimal immersion in $S^3$ that is not
totally geodesic, Simons' calculation~\cite{Simons} gives four
linearly independent eigenfunctions of the
Jacobi operator with eigenvalue $-2$. The simplicity of the lowest
eigenvalue then implies that the Morse index is at least five. In 1990 Urbano \cite{Urbano} characterized the equality case: the only
closed connected orientable minimal immersions of index at most
five are the equator and the Clifford torus, with indices one
and five, respectively. 
Perdomo~\cite[Theorem~4.3]{Perdomo} subsequently excluded index six
under antipodal symmetry. For closed orientable minimal surfaces
of genus $g\ge1$, Savo~\cite[Theorem~4.1]{Savo} proved
$I\ge4+\lceil g/2\rceil$. In view of the genus-zero and
genus-one classifications (see \cite{Almgren, Calabi} and \cite{Brendle} respectively), this confines a hypothetical embedded index-six
example to genera two, three, and four.

More recently, Kapouleas and Wiygul~\cite{KW} proved that Lawson's surface
$\xi_{g,1}$ (of genus $g$) has index $2g+3$ and nullity six for any $g\ge2$.
A short time ago, Karpukhin and Schulz~\cite[Theorem~1.3]{KS} were able to recover the index and
nullity of $\xi_{g,1}$ solely from topology and an index-two subgroup of
its symmetry group, while Kusner and Wang~\cite[Theorem~1.4]{KusW}
prove $\ind(\xi_{p,q})\ge2\max\{p,q\}+3$. 

Here we shall prove a refined bound that yields $I\ge7$ for every closed
orientable minimal immersion of genus at least two. In fact, this comes as a byproduct of an improvement of
Savo's aforementioned estimate. Throughout, we agree that the round three-sphere has sectional curvature one and we use the following notation. Let
$F:M\to S^3\subset\R^4$ be a smooth minimal immersion of a
closed connected orientable surface of genus $g$, and let
$h=F^*g_{S^3}$ be its induced metric. We write $K_h$ for the
Gauss curvature of $h$ and choose a smooth global unit normal
$\nu$ along $F$, so that $\nu(p)\in T_{F(p)}S^3$ is perpendicular
to $dF_p(T_pM)$ and has unit length. Let $A$ be the shape operator.
Our convention is
$\Delta=\diver\nabla$, so $-\Delta$ is nonnegative.
The Jacobi operator and the quadratic form associated to the second variation of the area functional are, respectively,
\begin{equation}\label{eq:J}
 J=-\Delta-|A|^2-2,
 \qquad
 Q(f)=\int_M\bigl(|df|^2-(|A|^2+2)f^2\bigr)\dd.
\end{equation}
Thus, for a normal variation with velocity $f\nu$, the second variation
of area is equal to $Q(f)$. We write $Q(f,g)$ for the associated symmetric
bilinear form. The index $I$ counts the negative eigenvalues of $J$ with
multiplicity, and the nullity is $N=\dim\ker J$.
An immersion is \emph{linearly full} if
$\operatorname{span}F(M)=\R^4$.

\begin{theorem}\label{thm:refined}
Let $b$ be the number of eigenvalues of $J$ in $(-2,0)$,
counted with multiplicity. If $g\ge2$, then
\begin{equation}\label{eq:refined}
 3I\ge2g+4b+N+9.
\end{equation}
In particular,
\begin{equation}\label{eq:genus-bound}
 I\ge5+\left\lceil\frac{2g}{3}\right\rceil.
\end{equation}
\end{theorem}

\begin{corollary}\label{cor:no-six}
There is no closed connected embedded minimal surface of Morse
index six in the round three-sphere.
\end{corollary}

\begin{remark}[Genus-two sharpness]
The general conclusion $I\ge7$ for surfaces of genus at least two is sharp since Lawson's
surface $\xi_{2,1}$ has Morse index equal to 7 (and nullity equal to $6$), by the index computation of Kapouleas and
Wiygul~\cite{KW}.
\end{remark}

The next classification question concerns index seven. Our bound
\eqref{eq:genus-bound}, together with the genus-zero and genus-one
classifications, restricts an embedded example to genus two or
three. We
propose the following:

\begin{conjecture}\label{conj:index-seven}
Up to ambient isometry, the only closed connected embedded minimal
surface of Morse index seven in the round three-sphere is
Lawson's surface $\xi_{2,1}$ of genus two.
\end{conjecture}

Let us now turn our attention to minimal immersions of nonorientable surfaces. For an orientable immersion, choosing a global unit normal
identifies normal sections with scalar functions on $M$.
For a nonorientable immersion, the normal line bundle is
nontrivial: normal sections instead correspond to odd scalar
functions on the orientation double cover. Simons' calculation~\cite{Simons} still gives the lower bound $I\ge4$. Shuli Chen~\cite[Theorem~3.6]{Chen} strengthened this (following Savo's approach), for
nonorientable minimal immersions, to
$I\ge4+\lceil g/6\rceil$, where $g$ denotes the genus of the
orientation double cover. Since $g\ge1$, this implies $I\ge5$.
It turns out that the same construction behind the proof of Theorem \ref{thm:refined}, applied
to the orientable double cover, gives
the following sharpened genus-dependent estimate.
Here $\pi:\widehat M\to M$ denotes the orientation double cover,
with nontrivial deck transformation $\tau$. The genus of a
nonorientable surface is understood to be the genus $g$ of
$\widehat M$. A smooth function
$f:\widehat M\to\R$ is called \emph{odd} if
$f(\tau(p))=-f(p)$ for every $p\in\widehat M$.

\begin{theorem}\label{thm:onesided}
Let $F:M\to S^3$ be a smooth minimal immersion of a closed connected
nonorientable surface of genus $g$. Set $h=F^*g_{S^3}$ and
$\widehat h=\pi^*h$, and let $\widehat A$ be the shape operator
of $\widehat F=F\circ\pi$. Let $I$ and $N$ be the
index and nullity on sections of the normal bundle of $F$, and let
$b_-$ count, with multiplicity, the eigenvalues in $(-2,0)$ of
$\widehat J=-\Delta_{\widehat h}-|\widehat A|_{\widehat h}^2-2$
on odd functions on $\widehat M$. Then
\begin{equation}\label{eq:onesided-refined}
 3I\ge g+4b_{-} + N+10.
\end{equation}
In particular,
\begin{equation}\label{eq:onesided-genus}
 I\ge5+\left\lceil\frac{g+1}{3}\right\rceil.
\end{equation}
\end{theorem}

\begin{corollary}\label{cor:onesided-six}
Every smooth minimal immersion $F:M\to S^3$ of a closed connected
nonorientable surface satisfies $\ind(F)\ge6$.
\end{corollary}

Let us now briefly comment on the proof of our main theorems. We use ambient components of one-forms as test functions for
the Jacobi form, following Ros~\cite{Ros}, Savo~\cite{Savo}, and
Ambrozio--Carlotto--Sharp~\cite{ACS}. However, the comparison underlying our argument is the Ejiri--Micallef
energy--area identity~\cite[Theorem~2.1 and equation~(2.27)]{EM},
applied to the Gauss map $\nu:M\to S^3$  (cf. computation
in the proof of~\cite[Proposition~1]{ACS}). More precisely, this relates the
energy variation of $\nu$ with velocity $fF+\alpha^\#$ to the
area variation of $F$ with normal velocity $f\nu$: see Proposition \ref{prop:square} and note that equality holds
precisely when $d\alpha=0$ and $\delta\alpha=2f$. We specialize that computation as follows: for each mean-zero function $f$ in the sum of the nonpositive
eigenspaces of $J$, we solve $\delta du_f=2f$ and take
$\alpha=du_f+\eta$, where $\eta$ is harmonic. That is the second ``canonical subspace'' of variations, the first being closely related to the variations in the coordinate directions that were singled out by Simons in \cite{Simons}. And, indeed, the main result in Section~\ref{sec:quadratic} implies 
that the ten-dimensional
space $\{B\nu:B\in\Sym(4,\R)\}$ is in direct sum with the preceding variations. The resulting space of $\R^4$-valued fields, obtained by projecting along a fixed ambient frame (as in Ros' work), has dimension at
least $I+N+2g+9$. Section~\ref{sec:proofs} shows that projecting each of their
four scalar components onto the sum of the eigenspaces of $J$
with eigenvalues at most $-2$ defines an injective linear map.
Its target consists of four copies of a scalar space of
dimension $I-b$. Comparing dimensions therefore gives
\[
 I+N+2g+9\le4(I-b),
\]
which is precisely Theorem~\ref{thm:refined}. Section~\ref{sec:proofs} also treats nonorientable immersions by carrying out
the construction with odd functions and one-forms on the
orientable double cover, although the case of genus 1 requires special care for there is no available counterpart to Brendle's characterization theorem.

Lastly, let us note that the same method also gives improved estimates for free boundary minimal
surfaces in the Euclidean unit ball $B^3$. For a compact connected
nonflat orientable immersion of genus $g$ with $b\ge1$ boundary
components, it yields $3I\ge2g+b+3$, improving the bound
$3I\ge2g+b-1$ of Ambrozio--Carlotto--Sharp~\cite{ACSfb}.
For a nonorientable immersion with $b$ boundary components and
orientation double cover of genus $g$, the corresponding estimate
is $3I\ge g+b+3$. One uses harmonic one-forms $\alpha$ satisfying
$d\alpha=\delta\alpha=0$ on $M$ and $\alpha(\mu)=0$
on $\partial M$, where $\mu$ denotes the outward unit
conormal to $\partial M$ in $M$, and adjoins the four fields $\nu,F_1\nu,F_2\nu,F_3\nu$;
in the nonorientable case, the harmonic forms are taken in the odd
subspace on the cover. 

\section{Geometric preliminaries and rotational Jacobi fields}\label{sec:preliminaries}

Throughout this article, we denote by $\nabla$ the Levi--Civita connection of $(M,h)$
and by $D$ ordinary differentiation of $\R^4$-valued fields
along $M$ (equivalently: $D$ is the flat connection on $M\times \R^4$ if we identify tangent fields with their images
under $dF$, which we shall do).

Set $S=|A|^2$. Minimality and the Gauss equation give
\begin{equation}\label{eq:geometry}
 A^2=\frac S2\,\mathrm{Id},\qquad K_h=1-\frac S2.
\end{equation}
For tangent vector fields $v,X$, we have
\begin{equation}\label{eq:Gauss}
 D_vF=v,\qquad D_v\nu=-Av,\qquad
 D_vX=\nabla_vX+\langle Av,X\rangle\nu-\langle v,X\rangle F.
\end{equation}
 Taking traces and using the Codazzi
equation yields
\begin{equation}\label{eq:coordinate}
 \Delta F=-2F,\qquad \Delta\nu=-S\nu.
\end{equation}
In particular, the four normal-coordinate functions satisfy
$J\nu_i=-2\nu_i$; see
Simons~\cite[Lemma~5.1.3]{Simons}. Here $\nu_i$ denotes the $i$-th component of $\nu$
in the standard basis of $\R^4$.

If the immersion is not totally geodesic, its Hopf differential is
holomorphic and not identically zero. This holds in particular
when $g\ge2$. Its zero set is finite and agrees with the zero set of
$A$. Consequently $S>0$ and $A$ is invertible on a dense open
subset of $M$.

\begin{lemma}\label{lem:full-killing}
If $g\ge2$, then $F$ is linearly full and $(M,h)$ has no
nonzero Killing vector field.
\end{lemma}
\begin{proof}
If $F$ were not linearly full, its image would lie in a great
two-sphere. The immersion $F:M\to S^2$ would then be a local
diffeomorphism. Compactness makes it a covering of $S^2$, so
connectedness and simple connectedness imply $M\cong S^2$,
a contradiction.

Let $Y$ be an intrinsic Killing field and $\varphi_t$ its flow.
Then $\varphi_t^*h=h$. Since $\varphi_0=\mathrm{id}$ and the
differentials $d\varphi_t$ are invertible and depend continuously
on $t$, the flow preserves orientation. Its maps are therefore
holomorphic for the complex structure determined by $h$ and the
orientation. Thus $Y^{1,0}$ is a holomorphic section of
$T^{1,0}M$.

A holomorphic line bundle $L$ of negative degree on a compact
Riemann surface has no nonzero holomorphic section
\cite[Chapter~V, Lemma~3.5, and Chapter~XI, \S2]{Miranda}. Indeed, a nonzero holomorphic section
$s$ has an effective zero divisor, and
\[
 \deg L=\deg\operatorname{div}(s)
       =\sum_{p\in M}\operatorname{ord}_p(s)\ge0.
\]
Since $\deg T^{1,0}M=2-2g<0$, it follows that $Y^{1,0}=0$
and hence $Y=0$.
\end{proof}

\begin{lemma}\label{lem:nullity}
Under the hypotheses of Theorem~\ref{thm:refined}, $N\ge6$.
\end{lemma}
\begin{proof}
For a constant skew-symmetric matrix $B\in\mathfrak{so}(4)$,
the ambient rotations $e^{sB}$ give a variation through minimal
immersions. Its normal component
\[
 \phi_B=\langle BF,\nu\rangle
\]
is therefore a scalar Jacobi field: $J\phi_B=0$.
We show that the map $B\mapsto\phi_B$ is injective.

If $\phi_B=0$, then $BF$ is tangent to the immersed surface,
since it is also perpendicular to $F$. Write $BF=dF(Y)$ for the
resulting smooth intrinsic vector field. The restriction of an
ambient Killing field that is tangent to an immersion is an
intrinsic Killing field. Explicitly, for tangent vectors $u,v$,
\[
 (\mathcal L_Yh)(u,v)
 =\langle Bu,v\rangle+\langle u,Bv\rangle=0.
\]
Lemma~\ref{lem:full-killing} gives $Y=0$, hence $BF=0$ everywhere.
Since $F$ is linearly full, this implies $B=0$.
The six-dimensional space $\mathfrak{so}(4)$ thus injects into
$\ker J$.
\end{proof}

\section{An energy--area second-variation identity}\label{sec:energy-area}

Let $\Omega^1(M)$ denote the space of smooth one-forms, and let
$\delta$ be the $L^2(M,h)$-adjoint of the exterior derivative $d$.
For $\alpha\in\Omega^1(M)$, write $\alpha^\#$ for its
metric-dual tangent field, identified with its image under $dF$.
Our sign convention gives
\[
 \delta\alpha=-\tr_h(\nabla\alpha),\qquad
 \delta du=-\Delta u.
\]
All tensor norms are taken with respect to $h$ and the Euclidean
inner product, and $L^2$ norms and inner products use $dA_h$.

We now regard the unit normal as the Gauss map
$\nu:(M,h)\to S^3\subset\R^4$; see~\cite[Section~10]{Lawson}. Since $d\nu=-dF\circ A$ and
$|d\nu|^2=S$, equation~\eqref{eq:coordinate} is the
sphere-valued harmonic map equation
$\Delta\nu+|d\nu|^2\nu=0$.
Away from the umbilics, $\nu$ is a minimal immersion, called
the polar immersion, with induced metric $(S/2)h$ and unit normal $F$.

Consider the energy functional
\[
 E_h(\psi)=\frac12\int_M|d\psi|^2\dd,
 \qquad \psi:(M,h)\to S^3.
\]
The domain metric $h$ is fixed. A variation is a smooth family
$\nu_t:M\to S^3$ with $\nu_0=\nu$; it need not arise from a
variation of the original immersion $F$. Its velocity belongs to
\[
 \mathcal V=\Gamma^\infty(\nu^*TS^3)
 =\{V\in C^\infty(M,\R^4):\langle V,\nu\rangle=0\}.
\]
Conversely, every $V\in\mathcal V$ is the velocity of such a
variation. Here the tangent space is taken at the value of the
Gauss map: explicitly,
\[
 T_{\nu(p)}S^3=\nu(p)^\perp
 =\R F(p)\oplus dF_p(T_pM),
\]
whereas $T_{F(p)}S^3=dF_p(T_pM)\oplus\R\nu(p)$.
Thus every $V\in\mathcal V$ has a unique decomposition
$V=fF+\alpha^\#$, with $f\in C^\infty(M)$ and
$\alpha\in\Omega^1(M)$. For $V,W\in\mathcal V$, set
\begin{equation}\label{eq:vector-form}
 \begin{split}
 \mathcal T(V,W)&=\int_M\bigl(\langle DV,DW\rangle-S\langle V,W\rangle\bigr)\dd,\\
 \mathcal T(V)&=\mathcal T(V,V).
 \end{split}
\end{equation}
Since $\nu$ is harmonic as a map into $S^3$, the second variation
of $E_h$ depends only on the velocity $V=\left.\partial_t\nu_t\right|_{t=0}$
and is
\begin{equation}\label{eq:energy-second}
 \left.\frac{d^2}{dt^2}\right|_{t=0}E_h(\nu_t)
 =\int_M\bigl(|DV|^2-S|V|^2\bigr)\dd
 =\mathcal T(V).
\end{equation}
Indeed, the sphere constraint gives
$\langle\nu,\left.\partial_t^2\nu_t\right|_{t=0}\rangle=-|V|^2$;
differentiating the energy twice and integrating by parts using
$\Delta\nu=-S\nu$ yields~\eqref{eq:energy-second}.
All variations here are constrained to $S^3$. For $V\in\mathcal V$, let $V_1,\ldots,V_4$ be its components
in a fixed orthonormal basis of $\R^4$. Since
$|DV|^2=\sum_{i=1}^4|dV_i|^2$ and
$|V|^2=\sum_{i=1}^4V_i^2$, we have
\[
 \mathcal T(V)
 =\sum_{i=1}^4\int_M\bigl(|dV_i|^2-SV_i^2\bigr)\dd
 =\sum_{i=1}^4\bigl(Q(V_i)+2\|V_i\|_{L^2}^2\bigr).
\]

The following identity is the Ejiri--Micallef energy--area
comparison~\cite[Theorem~2.1 and equation~(2.27)]{EM} applied to
the Gauss map, expressed in terms of the metric induced by the
original immersion. We give a direct proof in one-form notation.

\begin{proposition}\label{prop:square}
For every $f\in C^\infty(M)$ and $\alpha\in\Omega^1(M)$,
\begin{equation}\label{eq:square}
 \begin{split}
 \mathcal T(fF+\alpha^\#)
 =Q(f)+\|d\alpha\|_{L^2}^2+\|\delta\alpha-2f\|_{L^2}^2.
 \end{split}
\end{equation}
\end{proposition}
In other words, the preceding identity can be rephrased as
\[
\boxed{
\delta^2 E_h(\nu)[fF+\alpha^\#]
=
\delta^2\operatorname{Area}(F)[f\nu]
+\|d\alpha\|_{L^2}^2
+\|\delta\alpha-2f\|_{L^2}^2.
}
\]
Here the energy variation is taken among maps $\nu_t:M\to S^3$,
with the domain metric $h$ fixed and velocity $fF+\alpha^\#$.
The area variation is taken among immersions $F_t:M\to S^3$,
with normal velocity $f\nu$, and its second variation is $Q(f)$. Note that $\nu_t$ need not be the
Gauss map of $F_t$. The identity gives an energy--area comparison,
with equality precisely when $d\alpha=0$ and $\delta\alpha=2f$.

\begin{proof}
We first recall the tangential calculation in the proof of
\cite[Proposition~1]{ACS}, retaining the Hodge-energy terms
for an arbitrary one-form. Equations~\eqref{eq:Gauss} and
\eqref{eq:geometry} give
\[
 |D\alpha^\#|^2
 =|\nabla\alpha|^2+\left(1+\frac S2\right)|\alpha|^2.
\]
Since $K_h=1-S/2$, the integrated Bochner--Weitzenb\"ock
formula (see~\cite[Section~2.4, equation~(6)]{ACS}) yields
\begin{align}
 \mathcal T(\alpha^\#)
 &=\int_M\bigl(|\nabla\alpha|^2+K_h|\alpha|^2\bigr)\dd\notag\\
 &=\|d\alpha\|_{L^2}^2+\|\delta\alpha\|_{L^2}^2.
 \label{eq:tangential}
\end{align}
For the radial term, $D_v(fF)=df(v)F+fv$ gives
\[
 \mathcal T(fF)=Q(f)+4\|f\|_{L^2}^2.
\]
Finally, \eqref{eq:Gauss} and integration by parts give
\[
 \mathcal T(fF,\alpha^\#)
 =-\int_M\bigl(\langle df,\alpha\rangle+f\delta\alpha\bigr)\dd
 =-2\int_M f\delta\alpha\dd.
\]
Expanding $\mathcal T(fF+\alpha^\#)$ proves~\eqref{eq:square}.
\end{proof}

For a harmonic one-form $\eta$, \eqref{eq:tangential} gives
$\mathcal T(\eta^\#)=0$, cf. ~\cite[Proposition~1]{ACS}.

\section{Quadratic independence of the Gauss map coordinates}\label{sec:quadratic}

For a symmetric $4\times4$ constant matrix $B$ with real
coefficients (seen as an endomorphism of $\R^4$), let
$\tr_T B$ be its trace on the tangent two-plane
$T_p=dF_p(T_pM)$.

\begin{lemma}\label{lem:quadratic}
Assume that $F$ is not totally geodesic and that $h$ is not flat.
If $B\in\Sym(4,\R)$ satisfies
\[
 \langle B\nu,\nu\rangle=0\quad\hbox{on }M,
\]
then $B=0$. As a result, the ten functions
$\nu_i\nu_j$, $1\le i\le j\le4$, are linearly independent.
The hypotheses hold in particular for surfaces of genus $g\ge2$.
\end{lemma}
\begin{proof}
Let $v_1,v_2$ be a local $h$-orthonormal frame of $TM$,
identified with its image under $dF$. Applying $\Delta$
to the assumed relation and using $\Delta\nu=-S\nu$
(see \eqref{eq:coordinate}) and $D_{v_a}\nu=-Av_a$ (see \eqref{eq:Gauss}), we obtain
\[
 0
 =\Delta\langle B\nu,\nu\rangle
 =-2S\langle B\nu,\nu\rangle
   +2\sum_{a=1}^{2}\langle BAv_a,Av_a\rangle
 =S\,\tr_T B.
\]
The last equality follows from the assumed relation and
the identities $A^*=A$ and $A^2=(S/2)\mathrm{Id}$ (see \eqref{eq:geometry}).
Since $S>0$ on a dense set, continuity implies
$\tr_T B=0$ everywhere. 

Next: at each point of $M$, the vectors $F,v_1,v_2,\nu$
form an orthonormal basis of $\R^4$ and so computing the trace
of $B$ in this basis gives
\[
 \tr B
 =\langle BF,F\rangle
   +\sum_{a=1}^{2}\langle Bv_a,v_a\rangle
   +\langle B\nu,\nu\rangle
 =\langle BF,F\rangle+\tr_T B+\langle B\nu,\nu\rangle.
\]
The last two terms vanish, respectively by the preceding
calculation and by the assumed relation. Therefore
$\langle BF,F\rangle=\tr B$, thus
in particular, $\langle BF,F\rangle$ is constant on $M$. Using the product rule, the symmetry
and constancy of $B$, and the identities $\Delta F=-2F$
and $D_{v_a}F=v_a$, we obtain
\begin{multline*}
 0
 =\Delta\langle BF,F\rangle
 =2\langle B\Delta F,F\rangle
   +2\sum_{a=1}^{2}
     \langle BD_{v_a}F,D_{v_a}F\rangle \\
 =-4\langle BF,F\rangle
   +2\sum_{a=1}^{2}\langle Bv_a,v_a\rangle
 =-4\tr B+2\tr_T B
 =-4\tr B.   
\end{multline*}
Consequently $\tr B=0$, and hence
$\langle BF,F\rangle=0$ everywhere on $M$.

Differentiating the last relation shows that $BF$ is perpendicular
to $T$, as well as to $F$. Hence
\[
 BF=b\nu,\qquad b=\langle BF,\nu\rangle.
\]
On the dense set where $A$ is invertible, differentiating
$\langle B\nu,\nu\rangle=0$ shows that $B\nu$ is perpendicular
to $T$ and to $\nu$. Symmetry of $B$ then gives $B\nu=bF$, at which stage continuity extends this identity to all of $M$.

For any tangent vector $v$, differentiation of $BF=b\nu$ gives
\begin{equation}\label{eq:Bv}
 Bv=db(v)\nu-bAv.
\end{equation}
Taking its inner product with $\nu$, and using symmetry and
$B\nu=bF$, gives $db(v)=0$. Thus $b$ is constant. Differentiating
$B\nu=bF$ and using $D_v\nu=-Av$ now yields
\[
 BAv=-bv.
\]
On the other hand, \eqref{eq:Bv} applied to $Av$ gives
$BAv=-bA^2v$. Consequently
\[
 b(A^2-\mathrm{Id})=0.
\]
If $b\ne0$, then $A^2=\mathrm{Id}$ and $K_h=0$, contrary to
the nonflatness assumption. Therefore $b=0$. Equations
$BF=B\nu=0$ and \eqref{eq:Bv} show that $B$ annihilates an
orthonormal basis of $\R^4$ at any point, so $B=0$.

Finally, a linear relation among the ten products is exactly an
identity $\nu^TB\nu=0$ for a constant symmetric matrix, with the
off-diagonal coefficients divided by two. The asserted independence
follows.
\end{proof}

\section{Proofs of the improved index estimates}\label{sec:proofs}

\begin{proof}[Proof of Theorem~\ref{thm:refined}]
Let $\mathcal E$ be the sum of the nonpositive eigenspaces of $J$,
and let
\[
 \mathcal E_0=\left\{f\in\mathcal E:\int_M f\dd=0\right\}.
\]
Then $\dim\mathcal E=I+N$ and
$\dim\mathcal E_0\ge I+N-1$. Let
\[
 \HH^1(M)=\{\eta\in\Omega^1(M):d\eta=0,\ \delta\eta=0\}
\]
be the space of harmonic one-forms. By Hodge theory,
$\dim\HH^1(M)=2g$.
For $f\in\mathcal E_0$, let $u_f$ be the unique mean-zero solution
of
\[
 \delta du_f=2f,
\]
equivalently $\Delta u_f=-2f$. The assignment $f\mapsto u_f$ is
linear. The one-forms satisfying $d\alpha=0$ and
$\delta\alpha=2f$ are precisely $\alpha=du_f+\eta$ with
$\eta\in\HH^1(M)$. Thus the space
\[
 \mathcal W=\{fF+(du_f+\eta)^\#:
                    f\in\mathcal E_0,\ \eta\in\HH^1(M)\}
             \subset\mathcal V
\]
is parametrized injectively by $(f,\eta)$: the $F$-component
determines $f$, hence $u_f$, and then $\eta$. Therefore
\begin{equation}\label{eq:lowerdim}
 \dim\mathcal W=\dim\mathcal E_0+2g\ge I+N+2g-1.
\end{equation}
For $W=fF+(du_f+\eta)^\#\in\mathcal W$, both squares
in Proposition~\ref{prop:square} vanish, so
\begin{equation}\label{eq:W-form}
 \mathcal T(W)=Q(f)\le0.
\end{equation}

We now add ten fields built from the Gauss map coordinates. Define
\[
 \mathcal Z=\{B\nu:B\in\Sym(4,\R)\}.
\]
We take all ambient components below in a fixed orthonormal basis
of $\R^4$. Each coordinate of $Z=B\nu$ satisfies $JZ_i=-2Z_i$.
Integration by parts therefore gives
\begin{equation}\label{eq:Z-form}
 Q(Z_i,v)+2\langle Z_i,v\rangle_{L^2}=0
 \qquad (i=1,\ldots,4,\ v\in C^\infty(M)).
\end{equation}
Lemma~\ref{lem:quadratic} shows that $B\mapsto B\nu$ is
injective, so $\dim\mathcal Z=10$. It also shows that
$\mathcal Z\cap\mathcal V=\{0\}$: if $B\nu\in\mathcal V$,
then $\langle B\nu,\nu\rangle=0$, hence $B=0$.
Since $\mathcal W\subset\mathcal V$, their sum is direct and
\begin{equation}\label{eq:augmented}
 \dim(\mathcal W\oplus\mathcal Z)\ge I+N+2g+9.
\end{equation}
Using~\eqref{eq:Z-form} with $v=W_i$ and $v=Z_i$ for each $i$,
then the component expression for $\mathcal T$ and
\eqref{eq:W-form}, we obtain
\[
\sum_{i=1}^4
   \bigl(Q(W_i+Z_i)+2\|W_i+Z_i\|_{L^2}^2\bigr)
 =\sum_{i=1}^4
   \bigl(Q(W_i)+2\|W_i\|_{L^2}^2\bigr)
   =\mathcal T(W)=Q(f)\le0.
\]
Let $\mathcal E_{\le-2}$ be the sum of the eigenspaces of $J$
with eigenvalues at most $-2$; by the definition of $b$,
$\dim\mathcal E_{\le-2}=I-b$. Let $P_{\le-2}$ denote the $L^2$-orthogonal projection of scalar
functions onto $\mathcal E_{\le-2}$. Applying this projection
to each ambient component defines a linear map
\[
 \Pi:\mathcal W\oplus\mathcal Z
 \longrightarrow \mathcal E_{\le-2}^{\oplus4},
 \qquad
 \Pi(V)=\bigl(P_{\le-2}V_1,\ldots,P_{\le-2}V_4\bigr),
\]
where $\mathcal E_{\le-2}^{\oplus4}$ is the direct sum of four
copies of $\mathcal E_{\le-2}$.

We claim that $\Pi$ is injective. Indeed, suppose that
$\Pi(V)=0$. Each component $V_i$ is then orthogonal to every
eigenfunction of $J$ with eigenvalue at most $-2$. The spectral
decomposition of $J$ therefore gives
\[
 Q(V_i)+2\|V_i\|_{L^2}^2\ge0,
\]
with strict inequality whenever $V_i\ne0$. If $V\ne0$, at least
one component is nonzero, and summing yields
\[
 \sum_{i=1}^4\bigl(Q(V_i)+2\|V_i\|_{L^2}^2\bigr)>0,
\]
contradicting the preceding inequality on
$\mathcal W\oplus\mathcal Z$. Thus $\Pi$ is injective, and hence
\[
 \dim(\mathcal W\oplus\mathcal Z)
 \le \dim\bigl(\mathcal E_{\le-2}^{\oplus4}\bigr)
 =4(I-b).
\]

This is the componentwise dimension comparison
underlying Ros~\cite[Theorem~16]{Ros} and
Ambrozio--Carlotto--Sharp~\cite[Section~4, Propositions~3--4]{ACS}. Together with
\eqref{eq:augmented}, it yields
\[
 I+N+2g+9\le4(I-b),
\]
which is \eqref{eq:refined}. Finally, Lemma~\ref{lem:nullity}
gives $N\ge6$; together with $b\ge0$, this yields $3I\ge2g+15$.
Since $I$ is an integer,
\eqref{eq:genus-bound} follows.
\end{proof}

\begin{proof}[Proof of Theorem~\ref{thm:onesided}]
Write $\widehat F=F\circ\pi$, so that
$h=F^*g_{S^3}$ and $\widehat h=\pi^*h=\widehat F^*g_{S^3}$.
Choose the unit normal $\nu$ along $\widehat F$ determined by the
orientations of $\widehat M$ and $S^3$. The deck transformation
$\tau$ reverses the orientation of $\widehat M$ while leaving
$\widehat F$ unchanged, so it reverses this choice of normal:
\[
 \widehat F\circ\tau=\widehat F,\qquad \nu\circ\tau=-\nu.
\]
The normal fibers over $p$ and $\tau(p)$ are naturally the same
line in $\R^4$. Under this identification, the normal field $f\nu$
descends to $M$ precisely when
\[
 (f\nu)(\tau(p))=(f\nu)(p),
 \qquad\text{equivalently}\qquad
 -f(\tau(p))\nu(p)=f(p)\nu(p).
\]
Thus normal sections on $M$ correspond to odd scalar functions
on $\widehat M$. More precisely, if $NF$ denotes the normal line
bundle of $F$, every $V\in\Gamma^\infty(NF)$ lifts uniquely as
$\pi^*V=f\nu$ with $f\in C^\infty_{\mathrm{odd}}(\widehat M)$.
Define the scalar Jacobi form
$\widehat Q:C^\infty_{\mathrm{odd}}(\widehat M)\to\R$ by
\[
 \widehat Q(f)=\int_{\widehat M}
 \bigl(|df|_{\widehat h}^2-(|\widehat A|_{\widehat h}^2+2)f^2\bigr)
 \,dA_{\widehat h}.
\]
Since $\pi$ is a two-sheeted local isometry, the second
variation of area satisfies
\[
 \delta^2\operatorname{Area}_F(V,V)
 =\frac12\,\widehat Q(f),
 \qquad V\in\Gamma^\infty(NF),\quad \pi^*V=f\nu.
\]
The metric $\widehat h$ and the function $|\widehat A|_{\widehat h}^2$
are invariant under $\tau$, so $\widehat J$ commutes with pullback
by $\tau$. In particular,
\[
 f\circ\tau=-f
 \quad\Longrightarrow\quad
 (\widehat Jf)\circ\tau=-\widehat Jf.
\]
Consequently $I$ counts the negative eigenvalues of $\widehat J$
on odd functions, and $N$ is the dimension of its odd kernel;
see also~\cite[Section~4, equations~(4.1)--(4.2)]{UrbanoOneSided}.
All $L^2$ norms and inner products below use $dA_{\widehat h}$.
We use $\delta_{\widehat h}$ for the
codifferential and $\mathcal T$ for the energy second variation
of $\nu:(\widehat M,\widehat h)\to S^3$,
as defined in~\eqref{eq:vector-form}.

Before constructing the test fields, we show that $\widehat F$
is linearly full, that $g\ge1$, and that $\widehat h$ is nonflat.
Linear fullness excludes total geodesicity; together with
nonflatness, it verifies the hypotheses of Lemma~\ref{lem:quadratic},
which will supply the ten independent fields $B\nu$.
Linear fullness will also enter the nullity estimate at the end
of the proof.

If $\widehat F$ were not linearly full, its image, and hence the
image of $F$, would lie in a great two-sphere. The immersion
$F:M\to S^2$ would be a local diffeomorphism, and pulling back
the orientation of $S^2$ would orient $M$, a contradiction.
Thus $\widehat F$ is linearly full and is not totally geodesic.
If $g=0$, its Hopf differential would vanish, since a genus-zero
Riemann surface has no nonzero holomorphic quadratic differentials.
Minimality would then give $\widehat A=0$, a contradiction. Hence $g\ge1$.

Suppose next that $\widehat h$ is flat. By the classical
classification of flat minimal surfaces in $S^3$
(see~\cite[p.~340]{Lawson}), the image of $\widehat F$ is a
Clifford torus $T$. Viewed as a map $F:M\to T$, the immersion $F$ is itself a local
diffeomorphism, since both $M$ and $T$ have dimension two.
Pulling back the orientation of $T$ would again orient $M$,
a contradiction. Therefore $\widehat h$ is nonflat.

Let $\mathcal E_-$ be the sum of the nonpositive eigenspaces of
$\widehat J$ on odd functions, so $\dim\mathcal E_-=I+N$.
Every $f\in\mathcal E_-$ has mean zero with respect to
$dA_{\widehat h}$. Let $u_f$ be the unique mean-zero solution of
\[
 \delta_{\widehat h}du_f=2f.
\]
Since $\tau$ is an isometry, $-u_f\circ\tau$ solves the same
equation and also has mean zero. Uniqueness therefore gives
$u_f\circ\tau=-u_f$.

Denote by $\HH^1_-(\widehat M)$ the space of odd harmonic one-forms:
\[
 \HH^1_-(\widehat M)
 =\{\eta\in\Omega^1(\widehat M):
       d\eta=0,\ \delta_{\widehat h}\eta=0,\ \tau^*\eta=-\eta\}.
\]
The Hodge star for $\widehat h$ interchanges the even and odd
harmonic one-forms, because $\tau$ reverses orientation.
The full harmonic space has dimension $2g$, so
\[
 \dim\HH^1_-(\widehat M)=g.
\]

For a one-form $\alpha$ on $\widehat M$, it is convenient here to distinguish its intrinsic metric-dual field
$Y_\alpha=\alpha^{\#_{\widehat h}}$ from its ambient representative
$X_\alpha=d\widehat F(Y_\alpha)$. If $\alpha$ is odd, then for every
$v\in T_p\widehat M$,
\[
 \widehat h_{\tau(p)}
       \bigl(Y_\alpha(\tau(p)),d\tau_pv\bigr)
 =(\tau^*\alpha)_p(v)
 =-\widehat h_p\bigl(Y_\alpha(p),v\bigr).
\]
Since $\tau$ is an isometry, this implies
$Y_\alpha(\tau(p))=-d\tau_pY_\alpha(p)$.
Differentiating $\widehat F\circ\tau=\widehat F$ then gives
\[
 \begin{split}
 X_\alpha(\tau(p))
 &=d\widehat F_{\tau(p)}Y_\alpha(\tau(p))\\
 &=-d\widehat F_{\tau(p)}d\tau_pY_\alpha(p)
 =-d\widehat F_pY_\alpha(p)
 =-X_\alpha(p).
 \end{split}
\]
Thus all four scalar components of $X_\alpha$ are odd.

The construction in the previous proof now gives
\[
 \mathcal W_-=
 \{f\widehat F+X_{du_f+\eta}:
       f\in\mathcal E_-,\ \eta\in\HH^1_-(\widehat M)\},
 \qquad
 \dim\mathcal W_-=I+N+g.
\]
Every $W\in\mathcal W_-$ is an odd $\R^4$-valued field
perpendicular to $\nu$. Proposition~\ref{prop:square}, applied
on $(\widehat M,\widehat h)$, gives
\[
 \mathcal T(W)=\widehat Q(f)\le0.
\]
The space $\mathcal Z=\{B\nu:B\in\Sym(4,\R)\}$ also consists
of odd fields. By Lemma~\ref{lem:quadratic}, it has dimension ten
and intersects $\mathcal W_-$ trivially.
Each component of $Z\in\mathcal Z$ satisfies
$\widehat JZ_i=-2Z_i$. Applying~\eqref{eq:Z-form} on the cover,
and then the component expression for $\mathcal T$, yields
\[
 \begin{aligned}
 &\sum_{i=1}^4
       \bigl(\widehat Q(W_i+Z_i)+2\|W_i+Z_i\|_{L^2}^2\bigr)\\
 &\qquad=\sum_{i=1}^4
       \bigl(\widehat Q(W_i)+2\|W_i\|_{L^2}^2\bigr)
       =\mathcal T(W)=\widehat Q(f)\le0.
 \end{aligned}
\]
Let $\mathcal E_{-,\le-2}$ be the sum of the eigenspaces of
$\widehat J$ on odd functions with eigenvalues at most $-2$.
Its dimension is $I-b_-$. Projection of each scalar component
onto this space defines a linear map
\[
 \Pi_-:\mathcal W_-\oplus\mathcal Z
 \longrightarrow\mathcal E_{-,\le-2}^{\oplus4}.
\]
This map is injective: if $\Pi_-(V)=0$, the spectral decomposition
on odd functions gives
$\widehat Q(V_i)+2\|V_i\|_{L^2}^2\ge0$ for every $i$, with strict
inequality whenever $V_i\ne0$. A nonzero field in the kernel would
therefore contradict the preceding inequality.
The target has dimension $4(I-b_-)$, so comparing dimensions gives
\begin{equation}\label{eq:onesided-bound}
 I+N+g+10\le4(I-b_-),
\end{equation}
which rearranges to~\eqref{eq:onesided-refined}.

Next, in order to prove \eqref{eq:onesided-genus}, it remains to estimate the nullity $N$.
For $B\in\mathfrak{so}(4)$, consider the family of minimal immersions
$\widehat F_s=e^{sB}\widehat F$. Its velocity is $B\widehat F$,
and its normal component is
\[
 (B\widehat F)^\perp=\phi_B\nu,
 \qquad \phi_B=\langle B\widehat F,\nu\rangle.
\]
As in Lemma~\ref{lem:nullity}, differentiation of this family
gives $\widehat J\phi_B=0$. Moreover,
\[
 \phi_B(\tau(p))
 =\langle B\widehat F(p),-\nu(p)\rangle
 =-\phi_B(p).
\]
Thus the scalar Jacobi function $\phi_B$ is odd, whereas the
normal field $\phi_B\nu$ is $\tau$-invariant and descends to $M$.
Consequently the linear map
\[
 \Phi:\mathfrak{so}(4)\longrightarrow
       \ker\bigl(\widehat J|_{C^\infty_{\mathrm{odd}}(\widehat M)}\bigr),
 \qquad \Phi(B)=\phi_B,
\]
is well defined.

If $B\in\ker\Phi$, then $B\widehat F$ is perpendicular to both
$\widehat F$ and $\nu$, and is therefore tangent to the immersion.
There is a unique smooth tangent field $Y_B$ with
$d\widehat F(Y_B)=B\widehat F$. For tangent vectors $u,v$,
\[
 \begin{split}
 (\mathcal L_{Y_B}\widehat h)(u,v)
 &=\langle B\,d\widehat F(u),d\widehat F(v)\rangle
   +\langle d\widehat F(u),B\,d\widehat F(v)\rangle\\
 &=0,
 \end{split}
\]
since $B$ is skew-symmetric. Hence $Y_B$ is an intrinsic Killing
field of $(\widehat M, \widehat h)$. The map $B\mapsto Y_B$ on $\ker\Phi$ is injective:
if $Y_B=0$, then $B\widehat F=0$ everywhere, and linear fullness
forces $B=0$. Writing $\operatorname{Kill}(\widehat M,\widehat h)$
for the space of intrinsic Killing fields, the standard rank--nullity relation gives
\[
 N\ge\operatorname{rank}\Phi
 =6-\dim\ker\Phi
 \ge6-\dim\operatorname{Kill}(\widehat M,\widehat h).
\]

If $g\ge2$, Lemma~\ref{lem:full-killing}, applied to
$\widehat F:\widehat M\to S^3$, shows that this Killing
space is zero, so $N\ge6$.
If $g=1$, we show directly that its dimension is at most one.
Let $K=K_{\widehat h}$. Since $\widehat h$ is nonflat, $K$ is
not constant: otherwise Gauss--Bonnet would force $K=0$, which we have already ruled out.
Choose $p\in\widehat M$ with $dK_p\ne0$. Every Killing field $Y$
preserves curvature, hence $Y(K)=0$ and $Y(p)\in\ker dK_p$.

Evaluation at $p$ is injective on Killing fields. Indeed, if
$Y(p)=0$, differentiating $Y(K)=0$ in any direction
$v\in T_p\widehat M$ gives
\[
 dK_p\bigl((\widehat\nabla_vY)_p\bigr)=0,
\]
where $\widehat\nabla$ is the Levi--Civita connection of
$\widehat h$. The endomorphism
$v\mapsto(\widehat\nabla_vY)_p$ is skew-adjoint, and in dimension
two a nonzero skew-adjoint endomorphism is invertible.
Since $dK_p\ne0$, it follows that $(\widehat\nabla Y)_p=0$.
The restriction of a Killing field to a geodesic is a Jacobi field.
Its value and covariant derivative at $p$ now vanish, so the
Jacobi equation gives $Y=0$ along every geodesic from $p$,
and hence on the connected closed surface $\widehat M$.
Therefore evaluation embeds the Killing space into the
one-dimensional space $\ker dK_p$. This proves
$\dim\operatorname{Kill}(\widehat M,\widehat h)\le1$ and thus $N\ge5$
when $g=1$.

Finally, if $g\ge2$, the conditions $N\ge6$ and $b_-\ge0$
in~\eqref{eq:onesided-refined} give $3I\ge g+16$, so
\[
 I\ge\left\lceil\frac{g+16}{3}\right\rceil
   =5+\left\lceil\frac{g+1}{3}\right\rceil.
\]
If $g=1$, the estimate $N\ge5$ instead gives $3I\ge16$, hence
$I\ge6=5+\lceil(g+1)/3\rceil$. These two cases together prove
\eqref{eq:onesided-genus}.
\end{proof}

\section*{Acknowledgements}
This project has received funding from the European Research Council
(ERC) under the European Union's Horizon 2020 research and innovation
programme (grant agreement No.~947923).


\begin{thebibliography}{99}
\bibitem{Almgren}
F. J. Almgren, Jr., \emph{Some interior regularity theorems for
minimal surfaces and an extension of Bernstein's theorem},
Ann. of Math. (2) \textbf{84} (1966), no.~2, 277--292.

\bibitem{ACS}
L. Ambrozio, A. Carlotto, and B. Sharp,
\emph{Comparing the Morse index and the first Betti number of
minimal hypersurfaces}, J. Differential Geom. \textbf{108}
(2018), no.~3, 379--410.

\bibitem{ACSfb}
L. Ambrozio, A. Carlotto, and B. Sharp,
\emph{Index estimates for free boundary minimal hypersurfaces},
Math. Ann. \textbf{370} (2018), no.~3--4, 1063--1078.

\bibitem{Brendle}
S. Brendle, \emph{Embedded minimal tori in $S^3$ and the Lawson
conjecture}, Acta Math. \textbf{211} (2013), 177--190.

\bibitem{Calabi}
E. Calabi, \emph{Minimal immersions of surfaces in Euclidean spheres},
J. Differential Geom. \textbf{1} (1967), no.~1--2, 111--125.



\bibitem{Chen}
S. Chen, \emph{Index of minimal hypersurfaces in real projective
spaces}, Int. Math. Res. Not. IMRN (2024), no.~10, 8592--8608.

\bibitem{EM}
N. Ejiri and M. Micallef, \emph{Comparison between second variation
of area and second variation of energy of a minimal surface},
Adv. Calc. Var. \textbf{1} (2008), 223--239.


\bibitem{KW}
N. Kapouleas and D. Wiygul, \emph{The index and nullity of the Lawson
surfaces $\xi_{g,1}$}, Camb. J. Math. \textbf{8} (2020), 363--405.

\bibitem{KS}
M. Karpukhin and M. B. Schulz, \emph{Index and nullity of minimal
surfaces via representation theory}, preprint (2026).
\href{https://arxiv.org/abs/2610.03023}{arXiv:2610.03023}.

\bibitem{KusW}
R. Kusner and P. Wang, \emph{Willmore stability of minimal surfaces
in the $n$-sphere}, preprint (2026).
\href{https://arxiv.org/abs/2611.26238}{arXiv:2611.26238}.

\bibitem{Lawson}
H. B. Lawson, Jr.,
\emph{Complete minimal surfaces in $S^3$},
Ann. of Math. (2) \textbf{92} (1970), no.~3, 335--374.


\bibitem{Miranda}
R. Miranda, \emph{Algebraic Curves and Riemann Surfaces},
Graduate Studies in Mathematics, vol.~5,
American Mathematical Society, Providence, RI, 1995.


\bibitem{Perdomo}
O. Perdomo, \emph{Low index minimal hypersurfaces of spheres},
Asian J. Math. \textbf{5} (2001), no.~4, 741--750.

\bibitem{Ros}
A. Ros, \emph{One-sided complete stable minimal surfaces},
J. Differential Geom. \textbf{74} (2006), no.~1, 69--92.

\bibitem{Savo}
A. Savo, \emph{Index bounds for minimal hypersurfaces of the sphere},
Indiana Univ. Math. J. \textbf{59} (2010), 823--837.

\bibitem{Simons}
J. Simons, \emph{Minimal varieties in Riemannian manifolds},
Ann. of Math. (2) \textbf{88} (1968), no.~1, 62--105.

\bibitem{Urbano}
F. Urbano, \emph{Minimal surfaces with low index in the
three-dimensional sphere}, Proc. Amer. Math. Soc. \textbf{108}
(1990), 989--992.

\bibitem{UrbanoOneSided}
F. Urbano, \emph{Second variation of one-sided complete minimal
surfaces}, Rev. Mat. Iberoam. \textbf{29} (2013), no.~2, 479--494.
\end{thebibliography}
\end{document}